\documentclass{amsart}

\usepackage[english]{babel}

\usepackage{amsmath} 
\usepackage{amssymb}
\usepackage{amsthm} 

\usepackage[all,2cell,emtex]{xy}
\input xypic                                      

\theoremstyle{plain} 
\newtheorem*{thm}{Theorem}
\newtheorem*{cor}{Corollary}  

\theoremstyle{definition}                                         
\newtheorem*{defin}{Definition}

\theoremstyle{remark}                                             
\newtheorem*{rem}{Remarks}
\newtheorem*{ex}{Example}
\newtheorem*{exs}{Examples}

\newcommand{\ie}{\hbox{\emph{i.e.}}}

\newcommand{\C}{\mathcal{C}}
\newcommand{\D}{\mathcal{D}}
\newcommand{\E}{\mathcal{E}}
\newcommand{\R}[1]{\operatorname{\mathsf{R}\hbox{$#1$}}}
\newcommand{\sr}[1]{\operatorname{\mathsf{R}\hbox{$\scriptstyle#1$}}}
\newcommand{\RR}[1]{\operatorname{{\kern 1pt}\rlap{\vrule height -1.5pt depth 2pt width 5.0pt}\kern -1pt\mathsf{R}\hbox{$#1$}}}
\newcommand{\Rr}[1]{\operatorname{{\kern .7pt}\rlap{\vrule height -1.5pt depth 2pt width 3.5pt}\kern -.7pt\mathsf{R}\hbox{$\scriptstyle#1$}}}
\newcommand{\LL}[1]{\operatorname{{\kern .7pt}\rlap{\vrule height -1.5pt depth 2pt width 4.3pt}\kern -.6pt\mathsf{L}\hbox{$#1$}}}
\newcommand{\Ll}[1]{\operatorname{{\kern .7pt}\rlap{\vrule height -1.5pt depth 2pt width 3.pt}\kern -.6pt\mathsf{L}\hbox{$\scriptstyle #1$}}}
\newcommand{\al}{\alpha}
\newcommand{\bt}{\beta}
\newcommand{\g}{\gamma}
\newcommand{\dl}{\delta}
\newcommand{\e}{\varepsilon}
\newcommand{\h}{\eta}
\newcommand{\ee}{\rlap{\vrule height -2pt depth 2.5pt width 4.3pt}\varepsilon}
\newcommand{\hh}{\rlap{\vrule height -2pt depth 2.5pt width 4.3pt}\eta}
\newcommand{\ees}{\rlap{\vrule height -1pt depth 1.5pt width 2.3pt}\varepsilon}
\newcommand{\hhs}{\rlap{\vrule height -2pt depth 2.5pt width 2.9pt}\eta}

\newcommand{\toto}{{\hskip -2.5pt\xymatrixcolsep{1.3pc}\xymatrix{\ar[r]&}\hskip -2.5pt}}
\newcommand{\gmto}{{\hskip -2.5pt\xymatrixcolsep{1pc}\xymatrix{\ar[r]&}\hskip -2.5pt}}

\newcommand{\wt}[1]{\widetilde{ #1 }}

\begin{document}

\title{Quillen's adjunction theorem for derived functors, revisited}

\author{Georges MALTSINIOTIS}
\address{Institut de Math\'ematiques de Jussieu\\
Universit\'e Paris 7 Denis~Diderot\\ \hfill\break
\phantom{aa} Case Postale 7012\\
2, place Jussieu\\ 
F-75251 PARIS cedex 05\\
FRANCE}
\email{maltsin@math.jussieu.fr}
\urladdr{http://www.math.jussieu.fr/$\,\widetilde{ \ }$maltsin/}


\begin{abstract}
The aim of this paper is to present a very simple original, purely formal, proof of Quillen's adjunction theorem for derived functors~\cite{Qu}, and of some more recent variations and generalizations of this theorem~\cite{DHKS},~\cite{R-B}. This is obtained by proving an abstract adjunction theorem for ``absolute'' derived functors. In contrast with all known proofs, the explicit construction of the derived functors is not used.
\end{abstract}

\maketitle

We recall that for every category $\C$ and every class of arrows $W$ in $\C$, there exist a category $W^{-1}\C$, called
\emph{localization of $\C$ by $W$}, and a functor $P:\C\gmto W^{-1}\C$, the \emph{localization functor}, that carries arrows in $W$ into isomorphisms in $W^{-1}\C$, and universal for this property: for every category $\D$ and every functor $F:\C\gmto\D$ such that $F(W)$ is contained in the class of isomorphisms in $\D$, there is a unique functor $\wt{F}:W^{-1}\C\gmto\D$ such that $F=\wt{F}P$. 
$$
\xymatrixcolsep{2.6pc}
\xymatrix{
\C\ar[d]_{P}\ar[dr]^{F}
\\
W^{-1}\C\ar@{-->}[r]_(.55){\wt{F}}
&\D
}\quad
$$
The category $W^{-1}\C$ is obtained from $\C$ by formally inverting the arrows in $W$; it has the same objects as $\C$ and the morphisms are equivalence classes of ``composable zigzags'' of morphisms in $\C$, the arrows going in the ``wrong direction'' belonging to $W$ (see \cite{GZ}). The category $W^{-1}\C$ is not necessarily \emph{locally small} {\ie} the \emph{class} of arrows from an object to another is not in general a \emph{set}. This set theoretic problem can be ignored using Grothendieck's notion of universe~\cite{SGA4}. In this paper, in order to simplify, we introduce the following definition.

\begin{defin}
A \emph{localizer} is a pair $(\C,W)$, where $\C$ is a category and $W$ a class of arrows in $\C$ such that
the localized category $W^{-1}\C$ is locally small.
\end{defin}

\begin{exs}
If $\C$ is a \emph{small} category, then for every set $W$ of arrows in $\C$, the pair $(\C,W)$ is a localizer. If $\C$ is a Quillen model category and $W$ the class of weak equivalences in $C$, then $(\C,W)$ is a localizer~\cite[Ch.~I, 1.13, Th.~$1'$]{Qu}.
\end{exs}

Let $(\C,W)$ be a localizer, $P:\C\gmto W^{-1}\C$ the localization functor and \hbox{$F:\C\gmto\D$} an arbitrary functor. We recall that a \emph{right derived functor} of $F$ is a pair $(\R{F},\al)$, where $\R{F}:W^{-1}\C\gmto\D$ is a functor and $\al:F\gmto\R{F}\circ P$ a natural transformation,
$$
\UseTwocells
\xymatrixcolsep{1pc}
\xymatrixrowsep{.6pc}
\xymatrix{
\C\ar[ddd]_{P}\ar[dddrrr]^{F}
\\
\ddrrtwocell<\omit>{\al}
\\
\\
W^{-1}\C\ar[rrr]_{\sr{F}}
&&&\D
}
$$
satisfying the following universal property. For every functor $G:W^{-1}\C\gmto\D$, and every natural transformation $\g:F\gmto G\circ P$, there is a unique natural transformation $\dl:\R{F}\gmto G$ such that $\g=(\dl\star P)\,\al$. 
$$
\xymatrixcolsep{2.5pc}
\xymatrix{
F\ar[d]_{\al}\ar[dr]^{\g}
\\
\R{F}\circ P\ar[r]_(.48){\dl\,\star\,P}
&G\circ P
}
$$
This condition means exactly that the functor $\R{F}$ (together with the natural transformation $\al$) is a \emph{left} Kan extension of $F$ along the localization functor $P$. The pair $(\R{F},\al)$ is an \emph{absolute} right derived functor of $F$ if for every functor $H:\D\gmto\E$, the pair $(H\circ \R{F},H\star\al)$ is a right derived functor of $H\circ F$. 
An absolute right derived functor of $F$ is in particular a right derived functor of $F$ (take $H=1_\D$).
\smallbreak

\begin{ex}
If $\C$ is a Quillen model category and \hbox{$F:\C\gmto\D$} a functor that carries weak equivalences between fibrant objects  in $\C$ into isomorphisms in $\D$, then there exists an \emph{absolute} right derived functor $(\R{F},\al)$ of $F$. In order to prove this, we observe that Quillen's existence theorem for derived functors~\cite[Ch.~I, 4.2, Prop.~1]{Qu} implies that $F$ has a right derived functor $(\R{F},\al)$ constructed as follows. For every object $X$ in $\C$, choose a fibrant resolution $i_X:X\gmto X'$; then $\R{F}(X)=F(X')$ and \hbox{$\al_X=F(i_X):F(X)\gmto F(X')=\R{F}(X)$}. As for every functor $H:\D\gmto\E$, the functor $H\circ F$ carries weak equivalences between fibrant objects in $\C$ into isomorphisms in $\E$, the same construction gives a right derived functor of $H\circ F$ equal to $(H\circ \R{F},H\star\al)$,  which proves  the statement.
\end{ex}

Let $(\C,W)$ and $(\C',W')$ be two localizers, $P:\C\gmto W^{-1}\C$ and $P':\C'\gmto W'{}^{-1}\C'$ the localization functors, and $F:\C\gmto\C'$ an arbitrary functor. A \emph{total right derived functor} (resp.~\emph{absolute total right derived functor}) of $F$ is a pair $(\RR{F},\al)$, where \hbox{$\RR{F}:W^{-1}\C\gmto W'{}^{-1}\C'$} is a functor and $\al:P'\circ F\gmto\RR{F}\circ P$ a natural transformation , which is a right derived functor (resp.~an absolute right derived functor) of $P'\circ F$. 
$$
\UseTwocells
\xymatrixcolsep{2.5pc}
\xymatrixrowsep{2.4pc}
\xymatrix{
\C\ar[d]_{P}\ar[r]^{F}
\drtwocell<\omit>{\al}
&\C'\ar[d]^{P'}
\\
W^{-1}\C\ar[r]_{\Rr{F}}
&W'{}^{-1}\C'
}
$$
\smallbreak

\begin{ex}
If $\C$ and $\C'$ are two Quillen model categories and \hbox{$F:\C\gmto\C'$} a functor that carries weak equivalences between fibrant objects in $\C$ into weak equivalences in $\C'$, then there exists an \emph{absolute} total right derived functor of $F$. To prove this, we observe that if $P'$ is the localization functor from $\C$ to its localization  by the weak equivalences, then the functor $P'\circ F$ carries weak equivalences between fibrant objects in $\C$ into isomorphisms, and the statement is  a particular case of the previous example.
\end{ex}
\goodbreak

The notions of \emph{left derived functor}, of \emph{absolute left derived functor}, of \emph{total left derived functor}, and  of \emph{absolute total left derived functor} are defined in a dual way.
\goodbreak\goodbreak

\begin{thm}
Let $(\C,W)$ and $(\C',W')$ be two localizers, 
$$P:\C\toto W^{-1}\C\qquad\hbox{and}\qquad P':\C'\toto W'{}^{-1}\C'$$ 
the localization functors, 
$$F:\C\toto\C'\quad,\qquad G:\C'\toto\C\quad$$
a pair of adjoint functors, and
$$\e:F\circ G\toto1_{\C'}\quad,\qquad\h:1_\C\toto G\circ F\quad$$
the unit and counit of the adjunction. We suppose that the functor $F$ \emph{(resp.~$G$)} has an absolute total left \emph{(resp.}~right\emph{)} derived functor $(\LL{F},\al)$ \emph{(resp.~$(\RR{G},\bt)$)}. 
$$
\UseTwocells
\xymatrixcolsep{2.5pc}
\xymatrixrowsep{2.4pc}
\xymatrix{
\C\ar[d]_{P}\ar[r]^{F}
\drtwocell<\omit>{^\al}
&\C'\ar[d]^{P'}
\\
W^{-1}\C\ar[r]_{\Ll{F}}
&W'{}^{-1}\C'
}
\kern 50pt
\xymatrix{
\C'\ar[d]_{P'}\ar[r]^{G}
\drtwocell<\omit>{\bt}
&\C\ar[d]^{P}
\\
W'{}^{-1}\C'\ar[r]_{\Rr{G}}
&W^{-1}\C
}
$$
Then the pair of functors
$$\LL{F}:W^{-1}\C\toto W'{}^{-1}\C'\quad,\qquad\RR{G}:W'{}^{-1}\C'\toto W^{-1}\C$$
is a pair of adjoint functors, and we can choose  the unit and counit of the adjunction
$$\ee:\LL{F}\circ\RR{G}\toto1_{W'{}^{-1}\C'}\quad,\qquad\hh:1_{W^{-1}\C}\toto\RR{G}\circ\LL{F}\quad$$
in such a way that the two following squares commute.
$$
\xymatrixcolsep{1.3pc}
\xymatrix{
\LL{F}\circ P\circ G\ar[rr]^{\Ll{F}\star\,\bt}\ar[d]_{\al\,\star\,G}
&&\LL{F}\circ\RR{G}\circ P'\ar[d]^{\ees\,\star\,P'}
&&&\RR{G}\circ P'\circ F
&&\RR{G}\circ\LL{F}\circ P\ar[ll]_{\Rr{G}\star\,\al}
\\
P'\circ F\circ G\ar[rr]_{P'\star\,\e}
&&P'
&&&P\circ G\circ F\ar[u]^{\bt\,\star\,F}
&&P\ar[ll]^{P\,\star\,\h}\ar[u]_{\hhs\,\star\,P}
}
$$
\end{thm}

\begin{proof}
As $(\RR{G},\bt)$ (resp.~$(\LL{F},\al)$) is an absolute total right (resp.~left) derived functor of $G$ (resp.~of $F$), {\ie} an absolute right (resp.~left) derived functor of $P\circ G$ (resp.~of $P'\circ F$), the pair $(\LL{F}\circ\RR{G},\LL{F}\star\,\bt)$ (resp.~$(\RR{G}\circ\LL{F},\RR{G}\star\,\al)$) is a right (resp.~left) derived functor of $\LL{F}\circ P\circ G$ (resp.~of $\RR{G}\circ P'\circ F$). The universal property of right (resp.~left) derived functors says that for every functor \hbox{$H':W'{}^{-1}\C'\gmto W'{}^{-1}\C'$} \hbox{(resp.~$H:W^{-1}\C\gmto W^{-1}\C$)} and every natural transformation $\g':\LL{F}\circ P\circ G\gmto H'\circ P'$ (resp.~$\g:H\circ P\gmto\RR{G}\circ P'\circ F$), there exists a unique natural transformation 
$$\dl':\LL{F}\circ\RR{G}\toto H'\qquad\quad\hbox{( resp.\quad}\dl:H\toto\RR{G}\circ\LL{F}\,\ )\quad$$
such that
$$\g'=(\dl'\star P')(\LL{F}\star\,\bt)\qquad\quad\hbox{( resp.\quad}\g=(\RR{G}\star\,\al)(\dl\star P)\,\ )\quad.$$
This universal property applied to the functor $H'=1_{W'{}^{-1}\C'}$ (resp.~$H=1_{W^{-1}\C}$) and the natural transformation
$$
\xymatrixcolsep{2.3pc}
\phantom{\hbox{(resp.~\ }}\g'=(P'\star\e)(\al\,\star\,G):
\xymatrix{
\LL{F}\circ\, P\circ G\ar[r]^(.48){\al\,\star\,G}
&P'\circ F\circ G\ar[r]^(.63){P'\star\,\e}
&P'
}
\!\!=1_{W'{}^{-1}\C'}\circ P'
\quad$$
$$
\xymatrixcolsep{2.3pc}
\hbox{(resp.~\ }\g=(\bt\,\star\,F)(P\,\star\,\h):\,1_{W^{-1}\C}\circ P=\!
\xymatrix{
P\ar[r]^(.33){P\,\star\,\h}
&P\circ G\circ F\ar[r]^(.48){\bt\,\star\,F}
&\RR{G}\circ\, P'\circ F
}\ ),$$
$$
\UseTwocells
\xymatrixcolsep{2pc}
\xymatrixrowsep{.8pc}
\xymatrix{
\C'\ar[ddd]_{P'}\ar[dddrrr]^{\kern 5pt\Ll{F}\circ P\circ G}
\\
\ddrrtwocell<\omit>{\kern 15pt\Ll{F}\star\bt}
\\
\\
W'{}^{-1}\C'\ar[rrr]_{\Ll{F}\circ\Rr{G}}
&&&W'{}^{-1}\C'
}
\kern 30pt
\xymatrix{
\C'\ar[ddd]_{P'}\ar[dddrrr]^{\kern 5pt\Ll{F}\circ P\circ G}
\\
\ddrtwocell<\omit>{<-.4>\kern 35pt {}^{(P'\star\e)(\al\star G)}}
\\
\\
W'{}^{-1}\C'\ar[rrr]_{1_{W'{}^{-1}\C'}}
&&&W'{}^{-1}\C'
}
$$
$$
\UseTwocells
\xymatrixcolsep{2pc}
\xymatrixrowsep{.8pc}
\xymatrix{
\C\ar[ddd]_{P}\ar[dddrrr]^{\kern 5pt\Rr{G}\circ P'\circ F}
\\
\ddrrtwocell<\omit>{^\kern -15pt\Rr{G}\star\al}
\\
\\
W^{-1}\C\ar[rrr]_{\Rr{G}\circ\Ll{F}}
&&&W^{-1}\C
}
\kern 30pt
\xymatrix{
\C\ar[ddd]_{P}\ar[dddrrr]^{\kern 5pt\Rr{G}\circ P'\circ F}
\\
&\ddrtwocell<\omit>{^<1.2>\kern -30pt {}_{(\bt\star F)(P\star\h)}}
\\
\\
W^{-1}\C\ar[rrr]_{1_{W^{-1}\C}}
&&&W^{-1}\C
}
$$
implies the existence of a unique natural transformation
$$\ee:\LL{F}\circ\RR{G}\toto1_{W'{}^{-1}\C'}\qquad\hbox{(resp.\quad}\hh:1_{W^{-1}\C}\toto\RR{G}\circ\LL{F}\,\ )$$
such that
$$(P'\star\e)(\al\,\star\,G)=(\ee\,\star P')(\LL{F}\star\,\bt)\qquad
\hbox{(resp.\quad}
(\bt\,\star\,F)(P\star\,\h)=(\RR{G}\star\,\al)(\hh\,\star P)\,\ ).$$

It remains to prove that
$$(\RR{G}\star\,\ee)(\hh\,\star\RR{G})=1_{\Rr{G}}
\qquad\hbox{and}\qquad
(\ee\,\star\LL{F})(\LL{F}\star\,\hh)=1_{\Ll{F}}\quad.$$
$$
\xymatrixcolsep{3.5pc}
\xymatrixrowsep{1.5pc}
\xymatrix{
\RR{G}\ar[r]^(.4){\hhs\,\star\Rr{G}}
&\RR{G}\circ\LL{F}\circ\RR{G}\ar[r]^(.6){\Rr{G}\star\,\ees}
&\RR{G}
\\
\LL{F}\ar[r]^(.4){\Ll{F}\star\,\hhs}
&\LL{F}\circ\RR{G}\circ\LL{F}\ar[r]^(.6){\ees\,\star\Ll{F}}
&\LL{F}
}
\vrule depth 50pt width 0pt
$$
The uniqueness part of the universal property of derived functors, implies that it is enough to prove that
$$\Bigl[\bigl((\RR{G}\star\,\ee)(\hh\,\star\RR{G})\bigr)\star P'\Bigr]\bt=\bt
\qquad\hbox{and}\qquad
\al\Bigl[\bigl((\ee\,\star\LL{F})(\LL{F}\star\,\hh)\bigr)\star P\Bigr]=\al\quad.$$
Let us prove the first of these two equalities:
$$\begin{aligned}
\Bigl[\bigl(&(\RR{G}\star\,\ee)(\hh\,\star\RR{G})\bigr)\star P'\Bigr]\bt
 =(\RR{G}\star\,\ee\,\star P')(\hh\,\star\RR{G}\circ P')\bt=\cr
&=(\RR{G}\star\,\ee\,\star P')(\RR{G}\circ\LL{F}\star\,\bt)(\hh\star P\circ G)
 =\Bigl[\RR{G}\star\bigl((\ee\,\star P')(\LL{F}\star\,\bt)\bigr)\Bigr](\hh\star P\circ G)\cr
&=\Bigl[\RR{G}\star\bigl((P'\star\e)(\al\,\star\,G)\bigr)\Bigr](\hh\star P\circ G)
 =(\RR{G}\circ P'\star\e)(\RR{G}\star\,\al\,\star\,G)(\hh\star P\circ G)\cr
&=(\RR{G}\circ P'\star\e)\Bigl[\bigl((\RR{G}\star\,\al)(\hh\star P)\bigl)\star\,G\Bigr]
 =(\RR{G}\circ P'\star\e)\Bigl[\bigl((\bt\,\star\,F)(P\star\,\h)\bigl)\star\,G\Bigr]\cr
&=(\RR{G}\circ P'\star\e)(\bt\,\star\,F\circ G)(P\star\,\h\star\,G)
 =\bt(P\circ G\star\e)(P\star\,\h\star\,G)\cr
&=\bt\Bigl[(P\star\bigl((G\star\e)(\h\star\,G)\bigr)\Bigr]=\bt\quad.
\end{aligned}$$
The second equality is proved dually.
\end{proof}
\goodbreak\goodbreak

\begin{cor} \emph{\textbf{(Quillen's adjunction theorem for derived functors.)}}
Let $\C$ and $\C'$ be two Quillen model categories and
$$F:\C\toto\C'\quad,\qquad G:\C'\toto\C\quad$$
a pair of adjoint functors. We suppose that $F$ \emph{(resp.~$G$)} carries weak equivalences between cofibrant objects in $\C$ \emph{(resp.}~between fibrant objects in~$\C'$\emph{)} into weak equivalences in $\C'$ \emph{(resp.}~in $\C$\emph{)}. Then $F$ \emph{(resp.~$G$)} has a total left \emph{(resp.}~right\emph{)} derived functor $(\LL,\al)$ \emph{(resp.~$(\RR{G},\bt)$)}, and the functor $\LL{F}$ is a left adjoint of the functor~$\RR{G}$.
\end{cor}

\begin{proof}
The example preceding the theorem and its dual imply that these derived functors exist and are \emph{absolute} derived functors. Therefore the corollary  is a particular case of our theorem.
\end{proof}

\begin{rem}
In his book~\cite{Qu}, Quillen proves his adjunction theorem under the additional hypothesis that $F$ preserves  cofibrations and $G$ preserves fibrations [\hbox{\emph{loc.~cit.}} Ch.~1, 4.5 Th.~3], but our proof shows that these hypotheses are not necessary. In the more recent books of Hovey~\cite{Ho} and of Hirschhorn~\cite{Hi}, the Quillen's adjunction theorem is proved with even stronger hypotheses by requiring $(F,G)$ to be a \emph{Quillen adjunction} between \emph{closed} model categories, with \emph{functorial factorizations}. The variants of the adjunction theorem in a more general setting then Quillen model categories, obtained by Dwyer, Hirschhorn, Kan and Smith~\cite{DHKS}, or by Radulescu-Banu~\cite{R-B}, in his paper on \emph{Anderson-Brown-Cisinski model categories}, alias \emph{cat\'egories d\'erivables}~\cite{Ci}, are also implied by ``the adjunction theorem for absolute derived functors'', proved here. This theorem will be used, in a forthcoming paper with Bruno Kahn~\cite{BKGM}, to prove a generalization of the adjunction theorem of Radulescu-Banu, implying all known adjunction theorems for derived functors.
\end{rem}

\bigskip

\end{document}